\documentclass [a4 paper,12pt]{article}
\usepackage{mathrsfs}

\usepackage[usenames]{color}
\usepackage{bbm,amsmath,amssymb,amsfonts,amsthm,graphics,graphicx}

\newtheorem{lem}{Lemma}
\newtheorem{thm}[lem]{Theorem}

\def\pf{\noindent {\it Proof.} }

\title{The asymptotic number\\ of different rooted trees of a
tree\footnote{Supported by NSFC and the ``973" program.}}
\author{\small Xueliang Li, Yiyang Li, Yongtang Shi\\
\small Center for Combinatorics and  LPMC-TJKLC\\
\small Nankai University, Tianjin 300071, China\\
\small Email: lxl@nankai.edu.cn, liycldk@mail.nankai.edu.cn,
shi@nankai.edu.cn}
\date{ }

\begin{document}

\maketitle

\begin{abstract}
Let $\mathcal{T}_n$ be the set of trees with $n$ vertices. Suppose
that each tree in $\mathcal{T}_n$ is equally likely. We show that
the number of different rooted trees of a tree equals
$(\mu_r+o(1))n$ for almost every tree of $\mathcal{T}_n$, where
$\mu_r$ is a constant. As an application, we show that the number of
any given pattern in $\mathcal{T}_n$ is also asymptotically normally
distributed with mean $\sim \mu_M n$ and variance $\sim \sigma_M n$,
where $\mu_M, \sigma_M$ are some constants related to the given
pattern.
This solves an open question claimed in Kok's thesis.\\
[2mm] Keywords: tree; rooted tree; pattern; generating function; limiting distribution; automorphism\\
[2mm] AMS Subject Classification (2010): 05C05, 05C30, 05D40, 05A15,
05A16.
\end{abstract}

\section{Introduction}

A {\it pattern} $M$ is a given small tree. We say that $M$ {\it
occurs} in a tree $T$ if $M$ is a subtree of $T$ in the sense that
the degree of each internal  vertex (of degree more than one) of $M$
matches the degree of the corresponding vertex in $T$, while each
external vertex (of degree one) of $M$ matches a vertex of $T$ with
an arbitrary degree. Occasionally, we say that the pattern is in a
tree instead of that the pattern occurs in a tree for abbreviation.
Let $\mathcal{T}_n$ be the set of trees with $n$ vertices. If we use
$X_{n,M}(T)$ to denote the number of a given pattern $M$ in
$\mathcal{T}_n$, then $X_{n,M}(T)$ is a random variable with
probability
$$\text{Pr}(X_{n,M}=k)=\frac{t_{n,k}}{t_n},$$ where $t_{n,k}$ denotes the number
of such trees in $\mathcal{T}_n$ that the number of pattern $M$ in
each of the trees is $k$, and $t_n=| \mathcal{T}_n |$.

Moreover, let $\mathcal{R}_n$ be the set of rooted trees. We can
also consider the number of a given pattern in $\mathcal{R}_n$.
Denote $X_{n,M}(R)$ to be the random variable.

The main work of this paper is to show that some random variable
satisfies
$$\frac{Y_n-E(Y_n)}{\sqrt{Var(Y_n)}}\rightarrow_w \mathcal{N}(0,1),$$
where $\mathcal{N}(0,1)$ is the random variable with standard normal
distribution and $\rightarrow_w$ means weak convergence. We then
call this $Y_n$ {\it asymptotically normal}. Moreover, if
$$\frac{Y_n-\mu n}{\sqrt{\sigma n}}\rightarrow_w \mathcal{N}(0,1),$$
then $Y_n$ is asymptotically normal with mean $\sim \mu n$ and
variance $\sim \sigma n$. We refer to \cite{gkthesis} for details.

In fact, it was shown in \cite{cdkk} that in $\mathcal{R}_n$ the
number $X_{n,M}(R)$ of any given pattern is asymptotically normal
with mean $\sim \mu_M n$ and variance $\sim \sigma_M n$, where
$\mu_M$ and $\sigma_M$ are some constants corresponding to the given
pattern. But, for the set $\mathcal{T}_n$ there is no such a result
on normal distribution. In \cite{gk}, the authors proved that for
any given pattern in $\mathcal{T}_n$ the limiting distribution has a
density $(a+bt^2)e^{ct^2}$, where $a, b, c$ are some constants. The
mean and variance of the number of any given pattern are still
asymptotically $\mu_M n$ and $\sigma_M n$ where the constants are
the same as in $\mathcal{R}_n$. Clearly, if one shows that $b=0$,
then the distribution is normal. For some special patterns, such as
a star (or a node with a given degree) pattern \cite{dg}, a
double-star pattern \cite{LL}, and a path pattern \cite{gkthesis},
the corresponding limiting distributions were proved to be normal.
For some previous work we refer to Robinson and Schwenk \cite{rs}.
For more details, we refer to \cite{cdkk,gk, gkthesis, rs}. However,
for any given pattern Kok in his thesis \cite{gkthesis} claimed that
it seems much more difficult to demonstrate the normality. In this
paper, we will solve this problem from a new point of view which is
different from the existing ones.  We study the number of different
rooted trees of tree and get that for almost every tree of order $n$
the number of different rooted trees of the tree is $(\mu_r
+o(1))n$. Then, for any given pattern $M$, since we already knew
that the number of pattern $M$ in $\mathcal{R}_n$ is asymptotically
normal, as a consequence, in $\mathcal{T}_n$ the limiting
distribution for $M$ is also normal.

We organize this paper as follows. In Section 2, we will introduce
some basic knowledge that will be used in our proofs. In Section 3,
we will present the detailed proofs. We concentrate on the number of
different rooted trees of a tree. Section 4 is devoted to study the
limiting distribution for any given pattern.

\section{Preliminaries}

Analogous to patterns, for each tree $T$ we use $X_n(T)$ to denote
the number of different rooted trees of $T$. Clearly, $X_n(T)$ is
also a random variable in $\mathcal{T}_n$ and notice that
$X_n(T)\geq 1$. We introduce the following two functions:
$$t(x)=\sum_{n\geq1}t_nx^n,$$
$$t(x,u)=\sum_{n\geq1,k\geq 1}t_{n,k}x^nu^k,$$
where the coefficient $t_{n,k}$ denotes the number of such trees in
$\mathcal{T}_n$ that each of the trees has $k$ different rooted
trees. Clearly, $\sum_{k\geq 1} t_{n,k}= t_n$. We always assume that
every tree of $\mathcal{T}_n$ is equally likely. Then,
$\mbox{Pr}(X_n(T)=k)=\frac{t_{n,k}}{t_n} $.

Let $T_n$ be a tree in $\mathcal{T}_n$. An automorphism $\Phi$ of $T_n$ is
defined as
\begin{align}
\Phi: &\mbox{} v_i\rightarrow \Phi(v_i)\notag\\
&v_iv_j\rightarrow \Phi(v_i)\Phi(v_j)\notag
\end{align}
where $v_i$ and $v_j$ are any two vertices in $T_n$ and $v_iv_j$ is
an edge of $T_n$ joining vertices $v_i$ and $v_j$. We call that two
vertices $u$ and $v$ of $T_n$ are in the same vertex class if $u$
can be mapped to $v$ by some automorphism. Clearly, this sets up an
equivalence relation on the vertex set of $T_n$, and hence the
vertices in $T_n$ are partitioned into some classes. If we designate
every vertex in a same vertex class to be the root, we shall get the
same rooted tree. Exactly to say, suppose that $v_i$ can be mapped
to $v_j$ under an automorphism $\Phi$ of $T_n$, and $R_n^{v_i}$,
$R_n^{v_j}$ are the two rooted trees of $T_n$ which rooted at $v_i$
and $v_j$. One can easily illustrate that the automorphism $\Phi$ is
also an isomorphism that maps $R_n^{v_i}$ to $R_n^{v_j}$.

Hence, the number of different rooted trees of a tree is exactly the
number of vertex classes of the tree. Then, let $X_n(T)$ also
represent the number of vertex classes of the tree under
automorphisms. Therefore, we can similarly introduce the random
variable $X_n(R)$ of vertex classes on the space of rooted trees
$\mathcal{R}_n$.

If we consider $X_n(R)$ in $\mathcal{R}_n$, we also suppose that
each tree in $\mathcal{R}_n$ is equally likely. We can define
similar generating functions on $\mathcal{R}_n$, and let $r(x)$,
$r(x,u)$ be the related functions, respectively. One can see that
$r(x,1)=r(x)$. Suppose
$$r(x,u)=\sum_{n\geq1,k\geq1}r_{n,k}x^nu^k,$$
where $r_{n,k}$ is the number of rooted trees in $\mathcal{R}_n$
that have $k$ vertex classes. It follows that in $\mathcal{R}_n$,
$$\mbox{Pr}(X_n=k)=\frac{r_{n,k}}{r_n}, $$ where
$r_n=|\mathcal{R}_n|$.

We should notice that when we count the number of vertex classes of
a rooted tree, the root itself always forms a class with a single
vertex, since any automorphism on a rooted tree must map the root to
itself. That is a bit different from the case for non-rooted trees.

Furthermore, suppose that the convergence radius of $r(x)$ is $x_0$.
Otter \cite{ot} showed that $x_0$ satisfies that $r(x_0)=1$ and the
asymptotic expansion of $r(x)$ is
\begin{eqnarray}\label{rexpa0}r(x)=1-b_1(x_0-x)^{1/2}+b_2(x_0-x)
+b_3(x_0-x)^{3/2}+\cdots,\end{eqnarray} where $x_0\approx 0.3383219$
and $b_1\approx 2.6811266$. And, $t(x)$ has a similar expansion,
namely, \begin{equation}\label{t(x)}
t(x)=c_0+c_1(x_0-x)+c_2(x_0-x)^{3/2}+\cdots.\end{equation} Applying
the transfer theorems in \cite{fr} on Eqs.(\ref{rexpa0}) and
(\ref{t(x)}), we get that
\begin{align}\label{numtr}
t_n&\sim \frac{C x_0^{-n}}{n^{5/2}},\notag\\
r_n&\sim \frac{D x_0^{-n}}{n^{3/2}},\notag
\end{align}
where $C$ and $D$ are some constants. For this, we refer to
\cite{p,rs}. It has been showed that $C=0.5349\ldots$ and
$D=0.4399\ldots$. The book \cite{fr} gives us more details on the
transfer theorems.

In what follows, we first investigate $X_n$ in $\mathcal{R}_n$. To
start with, we need the following two lemmas. We refer to \cite{d,
gkthesis} for detailed information.

\begin{lem}\label{main1}
Suppose that $F(x,y,u)$ is an analytic function around $(x_0,y_0,1)$
such that $F(x_0,y_0,1)= y_0$, $F_{{y}}(x_0,y_0,1)=1$,
$F_{yy}(x_0,y_0,1)\neq 0$ and $F_x(x_0,y_0,1)\neq 0$. Then there
exist a neighborhood $U_0$ of $(x_0,1)$, a neighborhood $U_1$ of
$y_0$ and analytic functions $g(x,u)$, $h(x,u)$ and $f(u)$ which are
defined on $U_0$, such that the only solutions $y\in U_1$ with $y=
F(x,y,u)$ and $(x,u)\in U_0$ are given by
$y(x,u)=g(x,u)+h(x,u)\sqrt{1-\frac{x}{f(u)}}$. Furthermore,
$g(x_0,1)=y_0$ and $h(x_0,1)=\sqrt{\frac{
2f(1)F_x(x_0,y_0,1)}{F_{yy}(x_0,y_0,1)}}$.\qed
\end{lem}

\begin{lem}\label{main2}
Let $y(x,u)$ denote a function defined on a neighborhood $U$ of
$(x_0,1)$, and
$y(x,u)=F(x,y(x,u),u)=g(x,u)+h(x,u)\sqrt{1-\frac{x}{f(u)}}$. If
$y(x,1)$ is aperiodic, i.e., if from $y(x,1) = x^r\tilde{y}(x^d,1)$
with some power series $\tilde{y}(x,u)$ it follows that $d = 1$ and
all the Taylor coefficients of $F_y(x,y,u)$ are non-negative, then
there exists an $\eta> 0$ such that $y(x,u)$ can be analytically
continuous in
$$\widetilde{U}=\{(x,u): |x|<x_0 +\eta, |u|<1+\eta, \mathrm{arg}(x-f(u))\neq 0, x\neq
f(u)\}.$$ Moreover, if $y(x,u)=\sum y_{n,k}x^nu^k=\sum [y_n(u)]x^n$
and $y_{n,k}\geq 0$, then
$$y_n(u)=\frac{h(f(u),u)}{2\sqrt{\pi}n^{3/2}}f(u)^{-n}+O(\frac{f(u)^{-n-1}}{n^{5/2}}).$$
And if $h(f(1),1)\neq 0$, we get that $X_n$ is asymptotically normal
with mean $\sim \mu n$ and variance $\sim \sigma n$.\qed
\end{lem}

{\noindent {\bf Remark:}} In \cite{d} and \cite{gkthesis}, the
authors always assumed that all the Taylor coefficients of
$F(x,y,u)$ are non-negative and had the same conclusions. We can
extend this condition to the above one: all the Taylor coefficients
of $F_y(x,y,u)$ are non-negative. However, we can completely follow
the steps of the proofs in \cite{d} and \cite{gkthesis} without any
changes to show the above lemmas. Hence, the proof will not be
repeated here.

\section{The number of different rooted trees of a tree}

Now we concentrate on the number of vertex classes of a rooted tree.
Recalled that an automorphism of a rooted tree must map the root to
itself, which is a bit different from an automorphism of a
non-rooted tree, namely, the root always forms a vertex class with a
single vertex. We shall show that $X_n(R)$ is asymptotically normal
with mean $(\mu_r+o(1)) n$ and variance $(\sigma_r +o(1)) n$ in
$\mathcal{R}_n$.

In what follows, there appears an expression of the form $Z^*(S_n;
f(x,u))$ (or $Z(S_n;f(x))$), which is the substitution of the
counting series $f(x,u)$ (or $f(x)$) into the cycle index $Z(S_n)$
of the symmetric group $S_n$. This involves replacing each variable
$s_i$ in $Z(S_n)$ by $f(x^i,u)$ (or $f(x^i)$). For instance, if
$n=3$, then $Z(S_3)=(1/3!)(s_1^3+3s_1s_2+2s_3)$ and $Z(S_3;
f(x))=(1/3!)(f(x)^3+3f(x)f(x^2)+2f(x^3))$,
$Z^*(S_3;f(x,u))=(1/3!)(f(x,u)^3+3f(x,u)f(x^2,u)+2f(x^3,u))$. We
refer to \cite{hp} for details, where it was shown that
$$r(x)=x\cdot\sum_{n\geq 0}Z(S_n;r(x))=x\cdot
\text{e}^{\sum_{k=1}\frac{r(x^k)}{k}}.$$ The coefficient of $x^p$ in
$Z(S_n;r(x))$ is the number of rooted trees of order $p+1$ whose
roots have degree $n$. Multiplication of $Z(S_n;r(x))$ by $x$
corrects the power of $x$ so that $x^p$ in $xZ(S_n;r(x))$ is the
number of those trees with $p$ vertices. This expression
$Z(S_n;r(x))$ follows from the P\'{o}lya Enumeration Theorem
\cite{hp}.

Analogously, we take the same procedure for $r(x,u)$ in this paper.
But, here we should notice that if the same two copies of a rooted
tree with $k$ vertex classes connect to a root, then the number of
vertex classes of the new rooted tree is $k+1$, because there is
only one new class, i.e., the new root, which is different from the
procedure for calculating the number of a given pattern. Hence, we
use $r(x^k,u)$ to denote the generating function for $k$ copies of a
rooted tree. And we can get
\begin{equation}\label{rrr}r(x,u)=xu\cdot\mbox{
e}^{\sum_{k\geq1}\frac{1}{k}r(x^k,u)}+\psi(x,u),\end{equation} where
the modification term $\psi(x,u)$ is a  series with integral
coefficients.

For instance, suppose that the tree has a root of degree $2$. Then,
$xu\cdot Z^*(S_2, r(x,u))=xu\cdot \frac{1}{2}(r(x,u)^2+r(x^2,u))$.
We have $r_{n,k}$ choices to form a rooted tree with the same two
branches. In $r(x,u)^2$ the term $r_{n,k}x^{2n}u^{2k}$ denotes the
number of such rooted trees. We should notice that the power $2k$
must be corrected into $k$ which means that the number of vertex
classes is still $k$. But, in $r(x^2,u)$ the term $r_{n,k}x^{2n}u^k$
denotes the number of those rooted trees with two same branches.
Moreover, note that $r(x^2,u)$ and $r(x^2,u^2)$ are both with
integral Taylor coefficients. Hence, there must be a modification
term $\psi_2(x,u)$ such that $xu\cdot Z^*(S_2, r(x,u))+\psi_2(x,u)$
performs the generating function of the trees with roots of degree
2. Clearly, $\psi_2(x,u)$ is a series with integral coefficients. We
can see that for any number of branches, we have to modify the
function when counting the numbers corresponding to the cases that
some branches are the same. Hence, in general, expression
(\ref{rrr}) follows.

Let $y=r(x,u)$, and $F(x,y,u)=xu\cdot\mbox{
e}^{\sum_{k\geq1}\frac{1}{k}r(x^k,u)}+\psi(x,u).$ Here, the Taylor
coefficients of $x, y$ and $u$ may not be non-negative. Recall that
$r(x,1)=r(x)$, that is, $\psi(x,1)=0$. Recall also that there exists a
real number $x_0$ such that $r(x_0)=1$ and $x_0$ is the convergence
radius.

Then, we have $F(x_0,y(x_0,1),1)=1=y(x_0,1)=r(x_0,1)$ and
$$F_{y}(x,y,u)=xu\cdot\mbox{
e}^{\sum_{k\geq1}\frac{1}{k}r(x^k,u)},$$ which implies that
$F_y(x_0,y_0,1)=1$. Moreover, $F_{yy}(x_0,y_0,1)\neq 0$ and
$F_x(x_0,y_0,1)\neq 0$. That is, all the conditions of Lemma
\ref{main1} hold. Furthermore, we have that the Taylor coefficients
of $F_{y}(x,y,u)$ are non-negative, because in $F_y$ all
$r(x^k,u)$'s ($k\geq 2$) are within non-negative Taylor coefficients
and the expression of $F_y$ has an exponential form. Thus, by Lemma
\ref{main2} we have that the random variable $X_n(R)$ is
asymptotically normal with mean
$$E(X_n(R))\sim \mu_r n\mbox{  } (n\rightarrow \infty)$$ and variance
$$Var(X_n(R))\sim\sigma_r n\mbox{ } (n\rightarrow \infty),$$
where $\mu_r$ and $\sigma_r$ are some constants. Here, we just
concentrate on the rooted trees. Some researchers had considered the
number of vertex classes of other kind of trees, such as
phylogenetic trees \cite{bf}, and the conclusion also points to an
asymptotically normal distribution.

In this paper, we mainly focus on the overall property of a
probability space. Following the book \cite{BB}, we will say that
{\it almost every} (a.e.) graph in a graph space $\mathcal{G}_n$ has
a certain property $Q$ if the probability $\mbox{Pr}(Q)$ in
$\mathcal{G}_n$ converges to $1$ as $n$ tends to infinity.
Occasionally, we will say {\it almost all} instead of almost every.

From Chebyshev inequality
$$\mbox{Pr}\left[\big|X_n-E(X_n)\big|>
n^{3/4}\right]\leq \frac{Var X_n}{n^{3/2}}\rightarrow 0 \mbox{ as }
n\rightarrow \infty,$$ it follows that for almost all rooted trees,
$E(X_n)-n^{3/4}\leq X_n\leq E(X_n)+n^{3/4}$, namely,
$X_n=(1+o(1))E(X_n)$. We can get that
\begin{thm}\label{riso}
For almost all rooted trees in $\mathcal{R}_n$, the number of vertex
classes under automorphisms is $(\mu_r +o(1))n$.\qed
\end{thm}

Therefore, we can study the number of vertex classes in a tree. To
get the final result, we need another property as follows. We have
defined the number of vertex classes of a tree. And, we call a
vertex {\it fixed} if this single vertex forms a class.

\begin{lem}\label{TN2RN}
Almost every tree in $\mathcal{T}_n$ has more than
$\lfloor\frac{1}{24}n\rfloor$ fixed vertices.
\end{lem}
\pf We prove this result by contradiction. Suppose that
$\mathcal{T}'_n$ is a subset of $\mathcal{T}_n$ such that every tree
$T'_n$ in $\mathcal{T}'_n$ has at most $\lfloor\frac{1}{24}n\rfloor$
fixed vertices. We first show that these fixed vertices form a
subtree in $T'_n$. In fact, for any two fixed vertices $v_1$ and
$v_2$, they can only map to $v_1$ and $v_2$ among themselves,
respectively. Thus, any $(v_1,v_2)$-path maps to the
$(v_1,v_2)$-path under any automorphism. So, all vertices in a
$(v_1,v_2)$-path are fixed ones, that is, all the fixed vertices
form a connected subgraph of $T'_n$. Consequently, the fixed
vertices induce a subtree $T''_n$ of $T'_n$ and $|T''_n|\leq
\lfloor\frac{1}{24}n\rfloor$.

If $|T''_n|=0$, then the tree $T'_n$ has a symmetrical edge. The
structure of $T'_n$ is determined by one half of the vertices in
$T'_n$. Hence, the number of trees in $\mathcal{T}_n$ having a
symmetrical edge is at most $|\mathcal{T}_{\frac{n}{2}}|$, and
$\frac{|\mathcal{T}_{\frac{n}{2}}|}{|\mathcal{T}_{n}|}\rightarrow
0$, which completes the proof.

\begin{figure}[ht]
\begin{center}
\includegraphics[height=6cm]{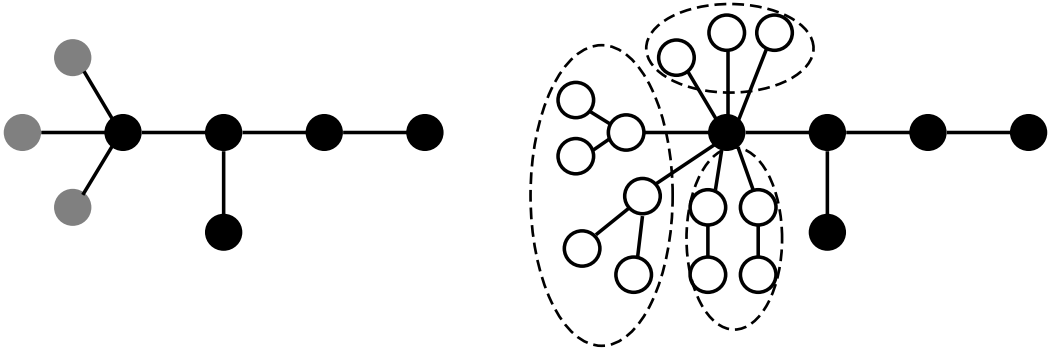}\\
\caption{An example of $A'$ and $A$.} \label{fig3}
\end{center}
\end{figure}

Then, we always suppose $|T''_n|>0$. Let $u$ be a vertex in $T''_n$.
Suppose that $H_u$ is a subtree of $T'_n$ attaching to $u$ such that
all the vertices in $H_u$ are not in $T''_n$. Suppose there are $m$
copies of $H_u$ after deleting $u$. We have $m\geq 2$; otherwise the
vertex in $H_u$ connecting to $u$ is also a fixed vertex, a
contradiction. If $m$ is even, we get rid of $m/2$ copies of $H_u$,
and if $m$ is odd, we get rid of $(m+1)/2$ copies of $H_u$. We
repeat this operation on all vertices in $T''_n$. At the end, this
produces a new tree $A$ with at most $\lfloor\frac{1}{2}(n+1/24\cdot
n)\rfloor$ vertices, and we denote the set of these new trees by
$\mathcal{A}_{\frac{25}{48}n}$. Moreover, if we replace these
$\lfloor \frac{m}{2}\rfloor$ copies of $H_u$ by a vertex, that is,
we add some vertices to $u$ and different kinds of $H_u$ correspond
to different vertices. Thus, we construct another tree $A'$. Observe
that $T''_n$ is a subtree of $A'$. We shall show that $A'$ has at
most $\lceil n/3\rceil$ vertices. Color the vertices in $A'$
corresponding to $T''_n$ by black and the others by gray. We already
knew that there are at most $\lfloor\frac{1}{24}n\rfloor$ black
vertices. Let $u$ be a black vertex. If a gray vertex connecting to
$u$ in $A'$ represents just only one vertex of $A$, then there is
only one such vertex connecting to $u$. One can see an example in
Figure \ref{fig3}. Hence, in $A'$ there are at most
$\lfloor\frac{1}{24} n\rfloor$ gray vertices representing the
subtrees of $A$ having a single vertex and at most
$\lfloor\frac{1}{2}\cdot \frac{1}{2} n\rfloor$ gray vertices
representing the other subtrees (or forests) having at least two
vertices. Consequently, we get that $A'$ has at most $\lfloor
n/3\rfloor$ vertices.

Moreover, we need the fact that the order of $\mathcal{T}_n$ is
asymptotically $\frac{C\cdot x_0^{-n}}{n^{5/2}}$. So, the number of
trees with at most $n$ vertices is asymptotically less than $2\cdot
C x_0^{-n}$. In the above, we have built a map from $\mathcal{T}'_n$
to $\mathcal{A}_{\frac{25}{48}n}$. Suppose $A$ is a tree in
$\mathcal{A}_{\frac{25}{48}n}$. Then $|A|$ is at most
$\lfloor\frac{25}{48}n\rfloor$. The ways of choosing $k$ vertices to
form a subtree of $A$ is at most ${\lceil\frac{25}{48}
n\rceil}\choose {k}$. We color these vertices in $A$ by black.
Notice that any tree in $\mathcal{T}''_n$ has at most
$\lfloor\frac{1}{24}n\rfloor$ vertices. Then, the number of all
subtrees within $k$ vertices in $A$ is less than
${{\lfloor\frac{25}{48} n\rfloor}\choose {k}} \leq
{{\lfloor\frac{25}{48} n\rfloor}\choose {\lfloor\frac{1}{24}}
n\rfloor}.$

 We select one subtree $T''_n$, and color the vertices black. Suppose that $A'$ is the corresponding
tree defined as above. For $u\in T''_n$, each gray vertex in $A'$
connecting to $u$ corresponds to a kind of subtree $H_u$. Moreover,
the number of $H_u$ can be odd or even in $T'_n$.  From the
structure of $A$, we reconstruct the tree $T'_n$ from $A$ by
deciding the number of $H_u$ to be odd or even. Since the number of
gray vertices is less than $|V(A')|\leq \lfloor n/3\rfloor$, we can
get that there exist at most $2^{\lfloor\frac{n}{3}\rfloor}$
different $T'_n$'s mapping to the same $A$.

Therefore, for trees in $\mathcal{T}''_n$ with $k$ vertices, at most
$2^{\lfloor\frac{n}{3}\rfloor}\cdot  2C\cdot x_0
^{-\frac{25}{48}n}\cdot {{\lfloor\frac{25}{48} n\rfloor}\choose {k}}
$ trees in $\mathcal{T}'_n$ map to them. Recall that each $T'_n$
corresponds to some $T''_n$. Then we have
\begin{align}
|\mathcal{T}'_n|&\leq \sum_{k=1}^{\lfloor1/24
n\rfloor}2^{\frac{n}{3}}\cdot 2C\cdot x_0 ^{-\frac{25}{48}n}\cdot
{{\lfloor\frac{25}{48} n\rfloor}\choose {k}}\notag
\\
&\leq \frac{1}{24}n 2^{\frac{n}{3}}\cdot 2C\cdot x_0 ^{-\frac{25}{48}n}\cdot {{\lfloor\frac{25}{48} n\rfloor}\choose {\lfloor\frac{1}{24} n\rfloor}}\notag\\
&=\frac{C}{12}n\cdot 2^{\frac{n}{3}}\cdot x_0
^{-\frac{25}{48}n}\cdot {{\lfloor\frac{25}{48} n\rfloor}\choose
{\lfloor\frac{1}{24} n\rfloor}}.\notag
\end{align}

By Stirling's approximation, i.e., $\frac{n!}{\sqrt{2\pi
n}{(\frac{n}{{\text{e}}}})^n}\rightarrow 1$ as $n\rightarrow
\infty$, we can get that when $n$ is large enough,
$${{\lfloor\frac{25}{48} n\rfloor}\choose {\lfloor\frac{1}{24} n\rfloor}}<
\frac{C_0}{\sqrt{n}}1.2^n,$$ where $C_0$ is a constant. Then,
$$|\mathcal{T}'_n|<C_1n^{1/2}\cdot
2^{\frac{n}{3}}x_0^{-\frac{25}{48}n}1.2^n,$$ where $C_1$ is some
real number for large $n$. It is known that
$|\mathcal{T}_n|\sim\frac{C\cdot x_0^{-n}}{n^{\frac{5}{2}}}$. Recall
that $x_0\approx 0.3383219$. Consequently,
$\frac{|\mathcal{T}'_n|}{|\mathcal{T}_n|}\rightarrow 0$.

Hence, in conclusion, we get that almost all trees do not belong to
$\mathcal{T}'_n$. The proof is thus complete.\qed
\\

Next, we proceed to estimate the number of different rooted trees of
a tree from Theorem \ref{riso} and Lemma \ref{TN2RN}. The following
theorem is established.

\begin{thm}\label{t2rt}
For almost all trees in $\mathcal{T}_n$, the number of different
rooted trees is $(\mu_r +o(1)) n$.
\end{thm}
\pf By Lemma \ref{TN2RN}, we know that almost every tree has at
least $\frac{1}{24}n$ fixed vertices, and denote these trees by
$\mathcal{T}^*_n$. Clearly, $\mathcal{T}^*_n \subseteq
\mathcal{T}_n$ and
$\frac{|\mathcal{T}^*_n|}{|\mathcal{T}_n|}\rightarrow 1$. Let $T$ be
a tree in $\mathcal{T}^*_n$. If we pick up one of the fixed vertices
to be the root, we can get a rooted tree having the same number of
vertex classes. There are at least $\frac{1}{24}n$ rooted trees in
which the roots of the rooted trees correspond to the fixed vertices
of $T$. And the number of vertex classes equals to that in $T$.
Hence, there are at least $|\mathcal{T}^*_n|\cdot \frac{1}{24} n$
rooted trees in $\mathcal{R}_n$ such that the roots are fixed
vertices in the associated tree. These rooted trees form a set
$\mathcal{R}^*_n$. Notice that $|\mathcal{R}_n|\sim \frac{D\cdot
x_0^{-n}}{n^\frac{3}{2}}$ and $|\mathcal{T}^*_n|\sim \frac{C\cdot
x_0^{-n}}{n^\frac{5}{2}}$. We get
$\frac{|\mathcal{R}^*_n|}{|\mathcal{R}_n|}\nrightarrow 0$. Combining
this with Theorem \ref{riso}, we have that the number of vertex
classes is $(\mu_r+o(1))n$ for almost all rooted trees in
$\mathcal{R}^*_n$.

According to whether the number of vertex classes is $(\mu_r+o(1))n$ or not, we depart
$\mathcal{R}^*_n$ into two parts $\mathcal{R}^*_{n,1}$ and
$\mathcal{R}^*_{n,2}$. There are at most
$\frac{|\mathcal{R}^*_{n,2}|}{\frac{1}{24}n}$ trees in
$\mathcal{T}^*_n$ corresponding to $\mathcal{R}^*_{n,2}$. Since
$|\mathcal{R}^*_{n,2}|=o(|\mathcal{R}^*_n|)=o(|\mathcal{R}_n|)$,
then
$\frac{|\mathcal{R}^*_{n,2}|}{\frac{1}{24}n}=o(|\mathcal{T}_n|)=o(|\mathcal{T}^*_n|)$.

Therefore, almost all trees in $\mathcal{T}^*_n$ correspond to the
rooted trees in $\mathcal{R}^*_{n,1}$. And recall that the root of
the tree in $\mathcal{R}^*_{n,1}$ is a fixed vertex. That is, almost
all trees in $\mathcal{T}^*_n$ also have $(\mu_r+o(1))n$ vertex
classes. Consequently, almost every tree in $\mathcal{T}_n$ has
$(\mu_r+o(1))n$ vertex classes. The proof is complete.\qed\\

From Theorem \ref{t2rt}, we have an intuitive grasp that the rooted
tree space is just the tree space with a scale $(\mu_r+o(1))n$. Not
rigorously to say, if we consider any special structure in trees,
the case that this structure will appear $(\mu_r+o(1))n$ times in
rooted trees is in a large probability, and the probabilities of
appearances in tree space and rooted tree space seem to be the same.
Moreover, by the asymptotical values of $|\mathcal{R}_n|$ and
$|\mathcal{T}_n|$, we can get that $\mu_r\approx 0.8210$.

\section{The distribution for any pattern in $\mathcal{T}_n$}

In this section, we shall focus on the distribution of the
occurrences $X_{n,M}(T)$ of a pattern $M$ on tree space
$\mathcal{T}_n$. It is known that the distribution of the
occurrences of a pattern in $\mathcal{R}_n$ is asymptotically
normal. We refer to \cite{gk} for this. We show that the
corresponding distribution in $\mathcal{T}_n$ is also asymptotically
normal. It has been shown that $X_{n,M}(T)$ has mean $(\mu_M+o(1))
n$ and variance $(\sigma_M +o(1))n$ and for almost every tree, and
the number of different rooted trees is $(\mu_r+o(1))n$. The
constants $\mu$ and $\sigma$ are the same as those for the case of
rooted trees, namely, $E(X_{n,M}(T))\sim\mu_M n\sim E(X_{n,M}(R))$
and $Var(X_{n,M}(T))\sim\sigma_M n \sim Var(X_{n,M}(R))$. Based on
these two results, we proceed to get our final result.

\begin{thm}
For any given pattern, the number of occurrences of the pattern in
trees is asymptotically normally distributed.
\end{thm}
\pf Recall that for each given pattern $M$, $E(X_{n,M}(T))\sim \mu_M
n$ and $Var(X_{n,M}(T))\sim \sigma n$, where $\mu_M$ and $\sigma_M$
are some constants. Let $\mathcal{T}_n^1$ be the subset of
$\mathcal{T}_n$ such that the number of occurrences $X_{n,M}(T)$
satisfies that $\frac{X_{n,M}(T)-\mu_M n}{\sqrt{\sigma_M n}}\leq t$,
where $t$ is some real number. Then, the probability
$$P\left(\frac{X_{n,M}(T)-\mu_M n}{\sqrt{\sigma_M n}}\leq t\right)=\frac{|\mathcal{T}_n^1|}{|\mathcal{T}_n|}, $$
where $\mathcal{T}_n^1$ is the subset of $\mathcal{T}_n$. For
$\mathcal{R}_n$,  we shall try to show that $$\lim_{n\rightarrow
\infty}P\left(\frac{X_{n,M}(R)-\mu_M n)}{\sqrt{\sigma_M n}}\leq
t\right) =\lim_{n\rightarrow \infty}P\left(\frac{X_{n,M}(T)-\mu_M
n}{\sqrt{\sigma_M n}}\leq t\right).$$

We knew that
\begin{align}
\lim_{n\rightarrow \infty}P\left(\frac{X_{n,M}(R)-\mu_M
n}{\sqrt{\sigma_M n}}\leq t\right)=N(0,1,t),\notag
\end{align}
where $N(0,1,t)$ denotes the probability value of the normal
distribution at $t$. Denote by
$\mathcal{R}_n^1$ the set of rooted trees satisfying
$\frac{X_{n,M}(R)-\mu_M n}{\sqrt{\sigma_M n}}\leq t$. The last equation holds from the fact that any
pattern in $\mathcal{R}_n$ is asymptotically normally distributed.

If $R$ is a rooted tree in $\mathcal{R}_n^1$ corresponding to $T\in
\mathcal{T}_n$, then $X_{n,M}(R)=X_{n,M}(T)$. So, a tree $T$ is in
$\mathcal{T}_n^1$ if and only if all the associated rooted trees are
in $\mathcal{R}_n^1$. We split $\mathcal{T}_n^1$ into two subsets,
$\mathcal{T'}_n^1$ and $\mathcal{T''}_n^1$, one is the collection of
trees corresponding to $(\mu_r +o(1))n$ rooted trees, and the other
is not, respectively. By Theorem \ref{t2rt}, the number of rooted
trees corresponding to $\mathcal{T'}_n^1$ is
$|\mathcal{T'}_n^1|\cdot(\mu(R) +o(1))n$, and
$|\mathcal{T''}_n^1|=o(|\mathcal{T}_n^1|)$, i.e., the number of
rooted trees associated with $\mathcal{T''}_n^1$ is at most
$o(|\mathcal{T}_n^1|)\cdot n$. Then, it follows that
$$|\mathcal{T'}_n^1|\cdot(\mu_r +o(1))n\leq |\mathcal{R}_n^1| \leq
|\mathcal{T'}_n^1|\cdot(\mu_r +o(1))n+o(|\mathcal{T}_n^1|)\cdot n.$$
Since $o(|\mathcal{T}_n^1|)\cdot n=o(|\mathcal{R}_n^1|)$ and
$\frac{|\mathcal{T'}_n^1|}{|\mathcal{T}_n^1|}\sim 1$, we have
$|\mathcal{R}_n^1|=(\mu_r +o(1))n\cdot |\mathcal{T}_n^1|$.
Therefore, we get that
\begin{align}
P\left(\frac{X_{n,M}(R)-\mu_M n}{\sqrt{\sigma_M n}}\leq t\right)&=\frac{\mathcal{R}^1_n}{\mathcal{R}_n}\notag\\
&\sim\frac{\mathcal{T}^1_n}{\mathcal{T}_n}\notag\\
&=P\left(\frac{X_{n,M}(T)-\mu_M n}{\sqrt{\sigma_M n}}\leq
t\right)\notag.
\end{align}
Consequently,
\begin{align}
\lim_{n\rightarrow \infty}P\left(\frac{X_{n,M}(T)-\mu_M
n}{\sqrt{\sigma_M n}}\leq t\right)
&=\lim_{n\rightarrow \infty}P\left(\frac{X_{n,M}(R)-\mu_M n}{\sqrt{\sigma_M n}}\leq t\right)\notag\\
&=N(0,1,t).\notag
\end{align}
Then the variable $X_{n,M}(T)$ is also asymptotically normal with
mean $\sim \mu_M n$ and variance $\sim \sigma_M n$. The proof is now complete. \qed\\

Now, we have established that for any pattern, the limiting
distribution of the number of occurrences in $\mathcal{T}_n$ is also
normal, which solves an open question claimed in \cite{gkthesis}.

\end{document}